\documentclass[12pt]{article}
\usepackage[T2A]{fontenc}
\usepackage[cp1251]{inputenc}
\usepackage[english]{babel}
\usepackage{amsthm, amssymb, amsmath
}
\begin{document}
\title{\normalsize 
\bf ON SYMMETRY OF EXTREMALS  \\
IN SEVERAL EMBEDDING THEOREMS}
\author { {\it E.V. Mukoseeva\footnote {The Chebyshev Laboratory, St. Petersburg State University} and 
A.I. Nazarov}\footnote {St. Petersburg Department of Steklov Mathematical Institute; 
Faculty of Mathematics and Mechanics, St. Petersburg State University }
} 
\date{}\maketitle

We consider the best constant problem in the following embedding theorem: 
\begin{equation}\label{problem}
\lambda(r,k,p,q)=\min\ \frac{\Vert f^{(r)}\Vert_{L_p\left[ -1,1\right] }}
{\Vert f^{(k)}\Vert _{L_q\left[-1,1\right] }}.
\end{equation}
Here $r,k\in\mathbb{Z}_+$, $r>k$, $1\le p,q\le\infty$, and minimum is taken over the set\footnote{The case of $p=1$ 
is special. In this case minimum should be taken over the set of $f$, which $(r-1)$-th derivative has bounded variation 
on $[-1,1]$. 
} ${\stackrel{\circ}{W}}\vphantom{W}_p^r(-1,1)$ that is
$$
\big\{f\in{\cal AC}^{r-1}[-1,1]\,\big|\,
f^{(r)}\in L_p(-1,1);\quad f^{(j)}(\pm1)=0,\ j=0,1,\dots,r-1\big\}.
$$

In the case of $r=1$, $k=0$ the problem (\ref{problem}) is well known. For $p=q=2$ it was solved by V.A. Steklov
\cite{St2}, for the arbitrary $p=q$ -- by V.I. Levin \cite{Lev} (see also \cite[Sect. 7.6]{HLP} for $p=q=2k$, $k\in\mathbb{N}$). 
Finally, E.Schmidt \cite{Sch} obtained the following result for arbitrary $p$ and $q$: 
$$
\lambda(1,0,p,q)=\frac{{\mathfrak F}\big( \frac 1q\big) {\mathfrak F}\big(\frac 1{p'}\big)}
{2^{\frac 1q + \frac 1{p'}}\,{\mathfrak F} \big( \frac 1q + \frac 1{p'} \big)},
$$
where ${\mathfrak F}(s) = \frac{\Gamma (s+1)}{s^s}$ and $p' = \frac{p}{p-1}$. Note
that the extremal in this problem is even function.\medskip

For $r=2$, $k=1$ the problem (\ref{problem}) is reduced to the best constant problem in the Poincare inequality, 
which was also solved by Steklov \cite{St} for $p=q=2$. However, the investigation
of the general case was completed only at the beginning of the XXI century and required efforts of many 
authors (\cite{DGS}, \cite{E}, \cite{BKN}, \cite{BeKa}, \cite{Kw}; the final result 
was obtained in \cite{N}). Namely, it turned out that for $q \leq 3 p$ the equality
$\lambda(2,1,p,q)=2\lambda(1,0,p,q)$ holds and the extremal is even function. However, for
$q >3 p$ we have $\lambda(2,1,p,q)<2\lambda(1,0,p,q)$ and the extremal is asymmetrical.\medskip

This result, as well as some calculations, leads to the following conjecture.\medskip

{\bf Conjecture}: {\it For $k\,\vdots\,2$ the extremal in the problem (\ref{problem}) is even function
for all admissible $r, p, q$. If $k\!\not\vdots\,2$ then for all admissible $r$ and $p$ 
there exists $\widehat q(r,k,p)>p$ such that the extremal is even for $q\le \widehat q$ and is asymmetrical for
$q>\widehat q$.}\medskip

For now, up to our knowledge, symmetry or asymmetry of the extremal is proved for the following parameters' values:

\medskip
\begin{center}
\begin{tabular} {|c||c|c|c|c||c|c|c|c|}
\hline \multicolumn{1}{|c||}{Article} & \multicolumn{4}{c||}{Symmetry} & \multicolumn{4}{c|}{Asymmetry} \\
\hline & $r$ & $ k$ & $p$ & $q$ & $r$ & $ k$ & $p$ & $q$ \\
\hline \cite{Sch} & $1$ & $0$ & $\forall$ & $\forall$ & \multicolumn{4}{c|}{ } \\
\hline \cite{BKN} & \multicolumn{4}{c||}{ } & $2$ & $1$ & $\forall$ & $>3p$ \\
\hline \cite{N} & $2$ & $1$ & $\forall$ & $\le 3p$ & \multicolumn{4}{c|}{ } \\
\hline \cite{NP}\footnotemark & $k+1$ & $\forall$ & $2$ & $2$ & \multicolumn{4}{c|}{ } \\
\hline \cite{Japans1} & $2,3$ & $0$ & $\forall$ & $\infty$ & \multicolumn{4}{c|}{ } \\
\hline \cite{Japans} & $\forall$ & $0$ & $2$ & $\infty$ & \multicolumn{4}{c|}{ } \\
\hline  \cite{Kal} & $\forall$ & $0,2$ & $2$ & $\infty$ & $\forall$ & $1$ & $2$ & $\infty$\\
\hline\end{tabular}
\end{center}
\footnotetext{see also \cite{J}.}
\medskip

In this paper we consider the case $p=2$, $q=\infty$. 
The main result is the following.\medskip

{\bf Theorem 1}. {\it Let $p=2$, $q=\infty$. 

\noindent {\bf 1}. If $k\!\not\vdots\,2$ then for all $r>k$ the extremal in the problem (\ref{problem}) 
is asymmetrical. 

\noindent {\bf 2}. If $k\,\vdots\,2$ then for all $r>k$ even function provides local minimum 
to the functional (\ref{problem}). }\medskip
 
We didn't manage to obtain complete solution for even $k$. The following theorem 
is developing results of \cite{Kal}.\medskip
 
{\bf Theorem 2}. {\it Let $p=2$, $q=\infty$. For $k=4,6$ and all $r>k$ the extremal in the problem
(\ref{problem}) is even. Furthermore,}
$$
\begin{aligned}
&\lambda(r,4,2,\infty)=\frac{1}{2^{r-2}(r-3)!}\sqrt{\frac{3(4r^2-24r+39)}{2(2r-9)}}; \\ 
&\lambda(r,6,2,\infty)=\frac{1}{2^{r-2}(r-4)!}\sqrt{\frac{192r^4-3456r^3+23372r^2-70240r+79065}{2(2r-13)}}.\\
\end{aligned}
$$
 
{\bf Proof of the theorem 1}.
Following \cite{Kal}, we introduce the function
 \begin{equation}\label{Kal}
A_{r,k}(x) = \max \{ |f^{(k)}(x)|: f\in{\stackrel{\circ}{W}}\vphantom{W}_2^r(-1,1),\Vert f^{(r)}\Vert_{L_2(-1,1)}\leq 1\}.
\end{equation}
Obviously, $\max\limits_{[-1,1]}A_{r,k}(x) =\lambda^{-1}(r,k,2,\infty)$.

We use the explicit formula, attained in \cite{Kal}: 
 \begin{equation}\label{Ark}
 A_{r, k}^2(x) = \big(Q_{r-k-1}(x)\big)^2\cdot\frac{1-x^2}{2(2r-2k-1)} - 
 \sum \limits_{n = r - k} ^ {r-1} {\big(Q_n^{(n+k-r)}(x)\big)^2\big(n + \frac 12\big)}, 
 \end{equation} 
where
$$
Q_n = \frac{1}{2^n n!} \cdot (1 - x^2)^n.
$$
Moreover, the function $f$ providing the maximum in (\ref{Kal}) is given by the following formula:
$$
f(t)=\sum\limits_{n\ge r} \big(n + \frac 12\big)\cdot Q_n^{(n+k-r)}(x)Q_n^{(n-r)}(t).
$$
It is easy to see that this function is symmetrical (even for $k$ even and odd for $k$ odd) if and only if $x=0$.

Thus, to prove {\bf Theorem 1} it is sufficient to prove the following lemma.\medskip 

{\bf Lemma}. {\it For odd $k$ the point $x=0$ provides the local minimum to the function $A_{r,k}(x)$. For even $k$ it provides the local maximum.}\medskip


Let us note that 
\begin{equation}\label{diff}
Q_n^{(s)}(0) = 
\begin{cases}
 0 & \hbox{for odd } s; \\
 \frac{(-1)^k s!}{2^n n!} C_n^k & \hbox{for } s = 2k.
\end{cases}
\end{equation}

Since the function $A_{r,k}$ is even we have $A_{r, k}'(0)=0$. Taking into account (\ref{diff}) we have
\begin{equation}\label{deriv}
\aligned
(A_{r, k}^2)''(0) = &-\frac{2}{(r-k-1)!^2 2^{2r-2k-1}} \\
  &- \sum \limits_{s = 0}^{k - 1} 2\big(r - k + s + \frac 12\big) \Big(\big(Q_{r - k + s}^{(s + 1)}(0)\big)^2  
+Q_{r - k + s}^{(s)}(0) Q_{r - k + s}^{(s + 2)}(0)\Big)  \\ 
= &-\frac{2}{(r-k-1)!^2 2^{2r-2k-1}} \\
  &- 2 \sum \limits_{t = 0}^{\lfloor\frac{k - 2}{2}\rfloor}
\Big(\frac{(2t+2)!C_{r - k + 2t + 1}^{t + 1}}{2^{r - k + 2t + 1} (r - k + 2t + 1)!} \Big)^2\big(r - k + 2t + \frac 32\big)\\ 
  & +2 \sum \limits_{t = 0}^{\lfloor\frac{k - 1}{2}\rfloor}
\frac{(2t)! (2t+2)! C_{r - k + 2t}^{t} C_{r - k + 2t}^{t + 1}}{2^{2(r - k + 2t)} (r - k + 2t)!^2}\big(r - k + 2t + \frac 12\big).
\endaligned
\end{equation}

Let $k$ be odd. In this case the second sum contains one additional term comparing to the first one. We separate the term corresponding to $t=0$ and get
\begin{equation*}
\aligned
(A_{r, k}^2)''(0) = \ & \frac{4(r - k)}{2^{2(r - k)} (r - k)!^2} \big(r - k + \frac 12\big)-\frac{2}{(r-k-1)!^2 2^{2r-2k-1}} \\ 
+ 2 \sum \limits_{t = 1}^{\frac{k - 1}{2}}\  & \bigg(\frac{(2t)! (2t+2)!C_{r - k + 2t}^{t} C_{r - k + 2t}^{t + 1}}
{2^{2(r - k + 2t)} (r - k + 2t)!^2} 
\big(r - k + 2t + \frac 12\big) \\ 
- & \Big(\frac{(2t)!C_{r - k + 2t -1}^{t}}{2^{r - k + 2t - 1} (r - k + 2t - 1)!} \Big)^2\big(r - k + 2t - \frac 12\big) \bigg). 
\endaligned
\end{equation*}
The expression in the first line is equal to $\frac{1}{2^{2r-2k-1} (r - k - 1)!^2 (r - k)}>0$. We denote it by $M$ and factor it out:
\begin{equation*}
\frac {(A_{r, k}^2)''(0)}M = 1 - 2 \sum \limits_{t = 1}^{\frac{k - 1}{2}}  
{\frac{(2t)!^2 (r-k-1)!^2(r-k)\big[(r-k)^2 + (r-k)(2t-1)-2t-\frac{1}{4}\big]}{2^{4t-1}t!^2(r-k+t-1)!^2(r-k+t)}}. 
\end{equation*}
The term in square brackets equals $(r-k+t-1)(r-k+t)-(t+\frac{1}{2})^2$, and we obtain
\begin{equation*}
\aligned
\frac {(A_{r, k}^2)''(0)}M = 1 - 2 \sum \limits_{t = 1}^{\frac{k - 1}{2}}  \ &
\Big(\frac{(2t)!^2 (r-k-1)!^2(r-k)(r-k+t-1)}{2^{4t-1}t!^2(r-k+t-1)!^2} \\
- \ & \frac{(2t)!^2 (r-k-1)!^2(r-k)(t+\frac{1}{2})^2}{2^{4t-1}t!^2(r-k+t-1)!^2(r-k+t)}\Big) \\ 
= 1 - 2 \sum \limits_{t = 1}^{\frac{k - 1}{2}} \ & \big(F(t)-F(t+1)\big)=1-2F(1)+2F\big(\frac {k+1}2\big),
\endaligned
\end{equation*}
where $F(t)=\frac{(2t)!^2 (r-k-1)!^2(r-k)(r-k+t-1)}{2^{4t-1}t!^2(r-k+t-1)!^2}$. Obviously, $2F(1)=1$, and therefore
\begin{equation*}
(A_{r, k}^2)''(0)=2MF\big(\frac {k+1}2\big) =
\frac{\big(k+1\big)!^2 \big(r-\frac{k+1}{2}\big)}{2^{2r-1}\big(\frac{k+1}{2}\big)!^2\big(r-\frac{k+1}{2}\big)!^2}>0,
\end{equation*}
which proves the first part of Lemma.\medskip

Now let $k=2\ell$ be even. Then the number of summands in both sums in (\ref{deriv}) equals $\ell-1$. 
Let us separate the term corresponding to $t=0$ from the second sum and add the term corresponding to $t=\ell$, which we subtract later. 
Then, similarly to the previous case, we get
\begin{equation*}
  (A_{r, k}^2)''(0) = M\cdot\Big(1-2\sum \limits_{t = 1}^{\ell}\big(F(t)-F(t+1)\big)\Big)-R=2MF\big(\ell+1\big)-R,
\end{equation*}
where $R = \frac{k!(k+2)!}{2^{2r-1}r!^2}C_r^\ell C_r^{\ell+1} \big(r+\frac 12\big)$.

After simplifying this expression we get
\begin{equation*}
\aligned
(A_{r, k}^2)''(0) & = \frac{(2\ell+2)!^2(r-\ell)}{2^{2r+1}(\ell+1)!^2(r-\ell)!^2}
- \frac{(2\ell)!(2\ell+2)!\big(r+\frac 12)}{2^{2r-1}\ell!(\ell+1)!(r-\ell)!(r-\ell-1)!} \\
& =-\frac{(2\ell)!(2\ell+2)!(r-\ell)}{2^{2r-1}\ell!(\ell+1)!(r-\ell)!(r-\ell-1)!}<0,
\endaligned
\end{equation*}
which proves the second part of the Lemma. Thus, Theorem 1 follows.\hfill$\square$\medskip

{\bf Proof of the theorem 2}. We use numerical-analytical method. Theoretically this scheme  
can be applied to any fixed $k$, but it requires more and more calculations when $k$ is increasing.\medskip

It is obvious from the formula (\ref{Ark}) that $A_{r,k}^2(x)=P_{r,k}(x^2)\cdot(1-x^2)^{2r-2k-1}$, where 
$P_{r,k}$ is a polynomial of degree $k$. Therefore
$$\frac{d\big[A_{r,k}^2(\sqrt{x})\big]}{dx}=P^{(1)}_{r,k}(x)\cdot(1-x)^{2r-2k-2},
$$
where $P^{(1)}_{r,k}$ is also a polynomial of degree $k$. It is easy to see that $P^{(1)}_{r,k}<0$ 
in a left semi-neighborhood of one. From the second statement of Theorem 1 we deduce that $P^{(1)}_{r,k}<0$
in a right semi-neighborhood of zero.

We construct a polynomial $\widetilde P_k$, such that all coefficients of the polynomial 
$\widetilde P_k(r\cdot)$ do not exceed the corresponding coefficients of $-P^{(1)}_{r,k}$. 
Thus, $-P^{(1)}_{r,k}(x)\ge\widetilde P_k(rx)$ for $x\ge0$. Then we show that 
the polynomial $\widetilde P_k$ is positive outside the interval $[c_1(k), c_2(k)]$. This means
that all roots of $P^{(1)}_{r,k}$ lye inside the interval $[\frac{c_1(k)}{r}, \frac{c_2(k)}{r}]$.\medskip

Now to proof the theorem it is sufficient to check that
\begin{equation}\label{<1}
\frac {P_{r,k}(x)}{P_{r,k}(0)}\cdot(1-x)^{2r-2k-1}<1,\qquad 
x\in\Big[\frac{c_1(k)}{r}, \frac{c_2(k)}{r}\Big].
\end{equation}

First, we prove (\ref{<1}) for $r$ big enough. To this end we rewrite $\frac {P_{r,k}(x)}{P_{r,k}(0)}$ as follows:
$$\frac {P_{r,k}(x)}{P_{r,k}(0)}=Q^+_{r,k}(x)-Q^-_{r,k}(x),
$$
where $Q^+_{r,k}$ is even polynomial and $Q^-_{r,k}$ is odd one. 

We construct polynomials $\widetilde Q^{\pm}_k$ with non-negative coefficients, such that for 
$r> r_0(k)$ all coefficients of the polynomial $\widetilde Q^+_k(r\cdot)$ are not less than the corresponding 
coefficients of $Q^+_{r,k}$ while the coefficients of $\widetilde Q^-_k(r\cdot)$ are not greater than the corresponding 
coefficients of $Q^-_{r,k}$. Then for $r> r_0(k)$ we have
\begin{equation*}
\aligned
\frac {P_{r,k}(x)}{P_{r,k}(0)}\cdot(1-x)^{2r-2k-1}
\leq(\widetilde Q^+_k(rx)-\widetilde Q^-_k(rx))\cdot(1-x)^{2(r-k)-1}\\
\leq(\widetilde Q^+_k(rx)-\widetilde Q^-_k(rx))\cdot\exp(-\alpha(k) rx),
\endaligned
\end{equation*}
where $\alpha(k)\le 2-\frac{2k+1}{r_0(k)}$.

Thus, the proof of (\ref{<1}) for $r> r_0(k)$ reduces to the proof of the following inequality:
$$
(\widetilde Q^+_k(y)-\widetilde Q^-_k(y))\cdot\exp(-\alpha(k)y)<1, \qquad y \in \big[c_1(k), c_2(k)\big]. 
$$
We prove this inequality by constructing suitable piecewise constant function $f_k$, which bounds 
the left-hand side of the inequality from above. To do so we note that the estimate
$$(\widetilde Q^+_k(y)-\widetilde Q^-_k(y))\cdot\exp(-\alpha(k)y)\le
(\widetilde Q^+_k(y_1)-\widetilde Q^-_k(y_0))\cdot\exp(-\alpha(k)y_0)
$$ 
holds for $c_1(k)\le y_0\le y\le y_1\le c_2(k)$, since the coefficients of polynomials $\widetilde Q^{\pm}_k$ are non-negative.

The obtained estimator $f_k$ was computed on a mesh fine enough. The inequality
$f_k<1$ proves (\ref{<1}) for $r> r_0(k)$.\medskip

We proceed similarly for $r\le r_0(k)$. Namely,
for every fixed $r\le r_0(k)$ we rewrite the polynomial in the left-hand side of (\ref{<1}) as follows:
$$\frac {P_{r,k}(x)}{P_{r,k}(0)}=R^+_{r,k}(x)-R^-_{r,k}(x),
$$
where $R^{\pm}_{r,k}$ are polynomials with non-negative coefficients.

We construct piecewise constant function $g_{r,k}$, which bounds the left-hand side of (\ref{<1}) from above.
Namely, for $\frac{c_1(k)}{r}\le x_0\le x\le x_1\le \frac{c_2(k)}{r}$ the following inequality holds: 
$$(R^+_{r,k}(x)-R^-_{r,k}(x))\cdot(1-x)^{2(r-k)-1}
\le(R^+_{r,k}(x_1)-R^-_{r,k}(x_0))\cdot(1-x_0)^{2(r-k)-1}.
$$

The obtained estimators $g_{r,k}$ were calculated on a mesh fine enough. The inequalities 
$g_{r,k}<1$ prove (\ref{<1}) for $r\le r_0(k)$, and the first statement of Theorem follows.

The values $\lambda(r,4,2,\infty)=\big(A_{r,4}(0)\big)^{-1}$ and $\lambda(r,6,2,\infty)=\big(A_{r,6}(0)\big)^{-1}$ are calculated
by the formulae (\ref{Ark}) and (\ref{diff}).\hfill$\square$\medskip

\subsection*{Appendix}

Here one can find the results of calculations, described in the proof of Theorem 2.

\paragraph{1. $k=4$.}

$$
\aligned
-P_{r,4}^{(1)}(x)&=(16r^4-96r^3+200r^2-168r+45)x^4+(-128r^3+656r^2-1056r+540)x^3\\
&+(312r^2-1224r+1134)x^2+(-240r+540)x+45;\\
\widetilde P_4(rx)&:=(3r^4)x^4+(-228r^3)x^3+(112r^2)x^2+(-350r)x+45.
\endaligned
$$
$$
c_1(4)=0.1;\qquad c_2(4)=76.
$$
Let us explain why the coefficients of $-P_{r,4}^{(1)}$ do not exceed the coefficients of $\widetilde P_4(r\cdot)$. 
Since $k=4$ implies $r\ge5$, it is sufficient to check that the difference of every corresponding coefficients pair  
is the polynomial in $r$ having a positive leading coefficient and no roots greater or equal than $5$. For instance, 
$$
(16r^4-96r^3+200r^2-168r+45)-3r^4 
=r^2(13r^2-96r+160)+(40r^2-168r+45).
$$ 
The roots of both quadratic polynomials in brackets are less than $5$. For other coefficients the argument is similar.

In the same manner we deduce that $\widetilde P_4$ is positive outside the interval $[c_1(4), c_2(4)]$: 
$$
3x^4-228x^3+112x^2-350x+45=(3x^2-228x+56)x^2+(56x^2-350x+45).
$$
Further, 
\begin{equation*}
\aligned 
&Q_{r,4}^+(x) = \Big(\frac{16}{9}r^4-\frac{128}{9}r^3+\frac{104}{3}r^2-32r+9\Big)x^4+\Big(\frac{56}{3}r^2-112r+126\Big)x^2+1;\\
&Q_{r,4}^-(x) = \Big(\frac{32}{3}r^3-\frac{688}{9}r^2+\frac{440}{3}r-84\Big)x^3+(8r-36)x;
\endaligned
\end{equation*}
\begin{equation*}
\aligned 
&\widetilde Q_4^+(rx) := \Big(\frac{17}{9}r^4\Big)x^4+\Big(\frac{57}{3}r^2\Big)x^2+1;\\
&\widetilde Q_4^-(rx) := \Big(\frac{26}{3}r^3\Big)x^3+(7r)x.
\endaligned
\end{equation*}
$$r_0(4)=50;\qquad \alpha(4)=1.8.
$$
Let us explain, why the coefficients of $Q_{r,4}^-$ are not less than the coefficients of $\widetilde Q_4^-(r\cdot)$ for 
$r>r_0(4)$. The difference of every corresponding coefficients pair is the polynomial in $r$, equals to the sum 
of binomials, positive for large $r$. For instance,
$$
\Big(\frac{32}{3}r^3-\frac{688}{9}r^2+\frac{440}{3}r-84\Big)-\frac{26}{3}r^3=
\Big(2r^3-\frac{688}{9}r^2\Big)+\Big(\frac{440}{3}r-84\Big).
$$
The value $r_0$ is chosen in the way that all the binomials are positive. One can show in the same way that coefficients of $\widetilde Q_4^+(r\cdot)$ 
are not less than the coefficients of $Q_{r,4}^+$.\medskip

The functions $f_4$ and $g_{r,4}$, $5\le r\le 50$, were computed on the following mesh:

\medskip
\begin{center}
\begin{tabular} {|l||c|}
\hline 
Function & Number of points\\
\hline  $f_4$ & $ 2^{11}\vphantom{\frac {1^1}{1^1}}$  \\
\hline  $g_{r,4}$, $5\le r\le 9$ & $2^{7}\vphantom{\frac {1^1}{1^1}}$  \\
\hline  $g_{r,4}$, $10\le r\le 21$ & $2^{8}\vphantom{\frac {1^1}{1^1}}$  \\
\hline  $g_{r,4}$, $22\le r\le 44$ & $2^{9}\vphantom{\frac {1^1}{1^1}}$  \\
\hline  $g_{r,4}$, $45\le r\le 50$ & $2^{10}\vphantom{\frac {1^1}{1^1}}$  \\
\hline\end{tabular}
\end{center}
\medskip
The computations were carried out with $17$ significant digits. This gives the estimate $1-f_4\ge 10^{-5}$, $1-g_{r,4}\ge 10^{-5}$.

\paragraph{2. $k=6$.}

\begin{equation*}
\aligned
-P_{r,6}^{(1)}(x)&=(64r^6-768r^5+3664r^4-8832r^3+11212r^2-6960r+1575)x^6\\
&+(-1152r^5+12768r^4-54528r^3+111792r^2-109560r+40950)x^5\\
&+(7440r^4-72960r^3+260760r^2-402240r+225225)x^4\\
&+(-21120r^3+170640r^2-449040r+386100)x^3\\
&+(26460r^2-156240r+225225)x^2+(-12600r+40950)x+1575; \\
\widetilde P(rx)&:=(4r^6)x^6+(-1200r^5)x^5+(1000r^4)x^4\\
&+(-22000r^3)x^3+(8000r^2)x^2+(-13000r)x+1575.
\endaligned 
\end{equation*}
$$c_1(6)=0.1;\qquad c_2(6)=300.
$$
Further,
\begin{equation*}
\aligned 
Q_{r,6}^+(x) &= \Big(\frac{64}{225}r^6-\frac{64}{15}r^5+\frac{5296}{225}r^4-\frac{1568}{25}r^3+\frac{19228}{225}r^2-\frac{836}{15}r+13\Big)x^6\\
        &+\Big(\frac{112}{5}r^4-288r^3+\frac{6104}{5}r^2-\frac{10584}{5}r+1287\Big)x^4\\
        &+(44r^2-396r+715)x^2+1;\\
Q_{r,6}^-(x) &= \Big(\frac{64}{15}r^5-\frac{1504}{25}r^4+\frac{22496}{45}r^3-\frac{17104}{25}r^2+\frac{10852}{15}r-286\Big)x^5\\
        &+\Big(\frac{736}{15}r^3-\frac{2736}{5}r^2+\frac{26216}{15}r-1716\Big)x^3+(12r-78)x.
\endaligned
\end{equation*}
\begin{equation*}
\aligned 
\widetilde Q^+(r\cdot) := &\Big(\frac{74}{225}r^6\Big)x^6+\Big(\frac{202}{9}r^4\Big)x^4+\Big(\frac{1982}{45}r^2\Big)x^2+1;\\
\widetilde Q^-(r\cdot) := &\Big(\frac{172}{45}r^5\Big)x^5+\Big(\frac{716}{15}r^3\Big)x^3+\Big(\frac{104}{9}r\Big)x.
\endaligned
\end{equation*}
$$r_0(6)=410;\qquad \alpha(6)=1.95.
$$
All the arguments for $k=4$ can be repeated word-by-word.
\medskip

The functions $f_6$ and $g_{r,6}$, $7\le r\le 410$, were computed on the following mesh:

\medskip
\begin{center}
\begin{tabular} {|l||c|}
\hline
Function & Number of points\\
\hline  $f_6$ & $2^{15}\vphantom{\frac {1^1}{1^1}} $  \\
\hline  $g_{r,6}$, $7\le r\le 42$ & $2^{11}\vphantom{\frac {1^1}{1^1}}$  \\
\hline  $g_{r,6}$, $43\le r\le 86$ & $2^{12}\vphantom{\frac {1^1}{1^1}}$  \\
\hline  $g_{r,6}$, $87\le r\le 173$ & $2^{13}\vphantom{\frac {1^1}{1^1}}$  \\
\hline  $g_{r,6}$, $174\le r\le 348$ & $2^{14}\vphantom{\frac {1^1}{1^1}}$  \\
\hline  $g_{r,6}$, $349\le r\le 410$ & $2^{15}\vphantom{\frac {1^1}{1^1}}$  \\
\hline\end{tabular}
\end{center}
\medskip

The computations were carried out with $17$ significant digits and gave the estimate $1-f_6\ge 10^{-5}$, $1-g_{r,6}\ge 10^{-5}$.\bigskip

The first author is supported by the Chebyshev Laboratory (St. Petersburg State University) 
under RF Government grant 11.G34.31.0026 and by JSC ``Gazprom Neft''. The second author is supported by RFBR grant 14-01-00534 and
by SPbSU grant 6.38.670.2013.

\end{document}